\documentclass{amsart}

\input{epsf.tex}

\theoremstyle{definition}

\theoremstyle{remark}

\newcommand {\bd}{\partial}
\newcommand{\leqslant}{\leq}
\newcommand{\geqslant}{\geq}
\newcommand{\G}{\gamma}
\newcommand {\da}{\delta}

\newcommand{\ra}{\rightarrow}

\numberwithin{equation}{section}

\begin{document}

\title{Asymptotic dimension and boundaries of hyperbolic spaces}

\author{Thanos Gentimis}
\address{Mathematics Department, University of Florida }

\email{thanos@ufl.edu}



\date{\today}


\begin{abstract}
We give an example of a visual Gromov-hyperbolic metric space $X$
with $asdim = 2$ and $dim(\partial X)=0$.
\end{abstract}

\maketitle

\
\section {Introduction}
The notion of asymptotic dimension of a metric space was
introduced by Gromov in \cite {GRO2}. It is a large scale analog
of topological dimension and it is invariant by quasi-isometries.
This notion has proved relevant in the context of Novikov's higher
signature conjecture. Yu \cite {Yu} has shown that groups of
finite asymptotic dimension satisfy Novikov's conjecture.
Dranishnikov (\cite {Dr }) has investigated further asymptotic
dimension generalizing several theorems from topological to
asymptotic dimension.

In this paper we are concerned with the relationship between
asymptotic dimension of a Gromov-hyperbolic space (see \cite {GRO1})
and the topological dimension of its boundary. Gromov in \cite
{GRO2}, sec. $1.E_1'$ sketches an argument that shows that complete
simply connected manifolds $M$ with pinched negative curvature have
asymptotic dimension equal to their dimension. He observes that the
same argument shows that $asdim(G)<\infty $ for $G$ a hyperbolic
group and asks whether such considerations lead further to the
inequality $asdim(G)\leqslant dim(\partial G)+1$.

Bonk and Schramm (\cite {Bo-S}) have shown that if $X$ is a
Gromov-hyperbolic space of bounded growth then $X$ embeds
quasi-isometrically to the hyperbolic $n$-space $\mathbb H $$^n$
for some $n$. It follows that $asdim (X)<\infty $ (see also \cite
{R} for a proof of this). If $K$ is any metric space one can
define (\cite {GRO1}, \cite {Bo-S}) a hyperbolic space $Con(K)$
with $\partial Con(K)=K$. If $X$ is a visual hyperbolic space then
$X$ is quasi-isometric to $Con(\partial X)$ (i.e. the boundary
'determines' the space). So it is natural to ask whether
$asdim(X)\leqslant \dim (\partial X)+1$ for visual hyperbolic
spaces in general. Besides the argument sketched in \cite {GRO2},
sec. $1.E_1'$ makes sense in this context too.

In this paper we give an example of a visual hyperbolic space $X$
 such that $asdim X = 2$ and $dim
\partial X = 0 $.
So the inequality $asdim(X)\leq dim (\partial X)+1$
 doesn't hold for this space.

We remark finally that Gromov's question for hyperbolic group was
settled in the affirmative recently by Buyalo and Lebedeva \cite
{B-L}.
\section{Preliminaries}

\textbf{Metric Spaces}. Let (X,d) be a metric space. The
\textit{diameter} of a set B is denoted by diam(B).
 A \textit{path} in X is a map $\G : I \ra X $ where I is an interval
in $\mathbb R$. A path $\G $ joins two points x and y in X if I [a,b] and $\G$(a) = x , $\G$(b) = y . The path $\G$ is called an
infinite ray starting from $x_0$ if I=[0,$\infty$) and $\G(0)=x_0$
. A geodesic, a geodesic ray or a geodesic segment in X is an
isometry $\G : I \ra X$ where I is $\mathbb R$ or $[0,\infty)$ or
a closed segment in $\mathbb R$. We use the term geodesic,
geodesic ray etc for the images of $\G$ without discrimination. On
a path connected space X given two points x,y we define the path
metric to be $\rho(x,y) = inf \{length(p)\}$ where the infimum is
taken over all paths $p$ that connect $x$ and $y$ (of course
$\rho(x,y)$ might be infinite). It is easy to see that inside a
ball B(x,n) of the hyperbolic plane or the euclidian plane the
path metric and the usual metric coincide. A metric space $(X,d)$
is called \textit{geodesic metric space} if $d=\rho $ (the path
metric is equal to the metric).

 \textbf{Hyperbolic Spaces}. Let (X,d) be a metric space .
Given three points x,y,z in X we define the \textit{Gromov
Product} of x and y with respect to the basepoint w to be
:$$(x|y)_w = \frac{1}{2} (d(x,w)+d(y,w)-d(x,y))$$  A space is said
to be \textit{$\da$- hyperbolic} if for all x,y,z,w in X we have:
$$(x|z)_w \geqslant min\{(x|y)_w,(y|z)_w\} - \da$$ A sequence of
points $\{x_i\}$ in X is said to converge at infinity if:
$$\lim_{i,j\ra \infty} (x_i|x_j)_w = \infty $$ Two sequences
$\{x_i\}$ and $\{y_i\}$ are equivalent if: $$\lim_{i,j\ra \infty}
(x_i|y_j) = \infty$$ This is an equivalence relation which does
not depend on the choice of w (easy to see). The boundary $\bd X$
of X is defined as the set of equivalence classes of sequences
converging at infinity. Two sequences are 'close' if
 $\liminf_{i,j\ra \infty}(x_i|y_j) $ is big. This defines a
 topology on the boundary.

 The boundary of every proper hyperbolic space is a
compact metric space.

 If $X$ is a geodesic hyperbolic metric space
and $x_0\in X$ then $\bd X$ can be defined as the set of geodesic
rays from $x_0$ where we define to rays to be equivalent if they are
contained in a finite Hausdorf neighborhood of each other. We equip
this with the compact open topology.

A metric $d$ on the boundary $\bd X$ of X is said to be
\textit{visual} if there are $x_0 \in X , a>1 $ and $ c_1,c_2 > 0$
such that $$c_1 a^{-(z,w)_{x_0}}\leqslant d(z,w) \leqslant
c_2a^{-(z,w)_{x_0}}$$ for every z,w in $\bd X$. The boundary of a
hyperbolic space always admits a visual metric (see \cite{GRO1}).

A hyperbolic space X is called \textit{visual} if for some $x_0
\in X$ there exists a $D>0$ such that for every $x \in X$ there
exists a geodesic ray $r$ from $x_0$ in $\bd X$ such that $d(x,r)
\leqslant D$ (see more on \cite {Bo-S}). It is easy to see that if
$X$ is visual with respect to a base point $x_0$ then it is visual
with respect to any other base point.
\newline

\textbf{Topological Dimension}. A covering $\{B_i\}$ has
\textit{multiplicity} $n$ if no more than $n + 1$ sets of the
covering have a non empty intersection. The \textit{mesh} of the
covering is the largest of the diameters of the $B_i$.

We will use in this paper the following definition of topological
dimension for compact metric spaces which is equivalent to the
other known definitions : A compact metric space has
\textit{dimension $\leqslant n$} if and only if it has coverings
of arbitrarily small mesh and order $\leqslant n$ . (see \cite
{H-W})\newline

\textbf{Asymptotic Dimension}.
 A metric space Y is said to be d - disconnected or that
it has dimension 0 on the d - scale if  $$ Y \,\, = \bigcup_{i\in I}
 B_i $$ such that: $sup \{ diam B_i , i \in I \}= D < \infty
 $, dist($B_i,B_j$) $\geqslant d$ $\forall  i \neq j$ where
dist($B_i,B_j$) = $\inf $ \{dist(a,b) a$\in B_i$ , b $\in B_j
\}$\newline
 \textit{(Asymptotic Dimension 1)}. We say that a space X has asymptotic
dimension n if n is the minimal number such that for every $d > 0$
 we have : $X =\bigcup X_k$ for k = 1,2, ... n and all $X_k$ are d-disconnected. We then write asdim = n

 We say that a covering  $\{B_i\}$  has \textit{d - multiplicity} ,
k if and only if every d - ball in X meets no more than k sets $B_i$
of the covering.A covering has \textit{multiplicity} n if no more
than n + 1 sets of the covering have one a non empty intersection. A
covering $\{ B_i \} \,i \in I$ is \textit{D - bounded} if diam
$(B_i) \leqslant D \,\forall \,i\in I$
\newline \textit{(Asymptotic Dimension 2)}. We say that a space $X$ has asdim
= n if n is the minimal number such that $\forall \,d > 0$ there
exists a covering of X of uniformly  D - bounded sets $B_i$ such
that d - multiplicity of the covering $\leqslant n+1$. The two
definitions are equivalent. (see \cite {GRO1})\newline

 \textbf{The Hyperbolic Plane}. The hyperbolic plane
$\mathbb H$$^2$ is a visual hyperbolic space of bounded geometry.
It is easy to see that $asdim\mathbb H$$^2 = 2$ (see \cite
{GRO2}).
We will use  the
standard model of the hyperbolic plane given by the interior of a
disk in $\mathbb R
$$^2$.

\section{Constructing The "COMB" Space}
Let $\mathbb {H}$$^2$ be the hyperbolic plane and let
$a_1,a_2,....$ be  geodesic rays starting from a point $x_0$ and
extending to infinity such that the angle between $a_n,a_{n+1}$ is
$\frac{\pi}{2^n}$.

Let $S(a_n,a_{n+1})$ be the sector defined by the rays
$a_n,a_{n+1}$. In other words $S(a_n,a_{n+1})$ is the convex
closure of $a_n,a_{n+1}$.

 Since geodesics diverge in $\mathbb {H}$$^2$ there
is an $x\in S(a_n,a_{n+1})$ such that the ball of radius $n$ and
center $x$, $B(x,n)$ is contained in $S(a_n,a_{n+1})$. Let $N_n$
be such that $B(x,n)\subset B(x_0,N_n)$. Let
$$S(a_n,a_{n+1},N_n)=S(a_n,a_{n+1})\cap B(x_0,N_n)$$

Let's call $K_n$ the upper arc of $S(a_n,a_{n+1},N_n)$, i.e.
$$K_n=S(a_n,a_{n+1},N_n)\cap \partial B(x_0,N_n)$$

 We subdivide
$K_n$ into small pieces of length between $1/2$ and 1 marking the
vertices.
Then we consider the geodesic rays starting from $x_0$
to every vertex we defined and we extend them to infinity.

So we arrive at the "COMB" space which is the union of all the
$S(a_n,a_{n+1},N_n)$ together with these rays and looks like this:

\begin{figure}[htb]
\epsfxsize 10.0cm \epsfysize 8.0cm \epsfbox{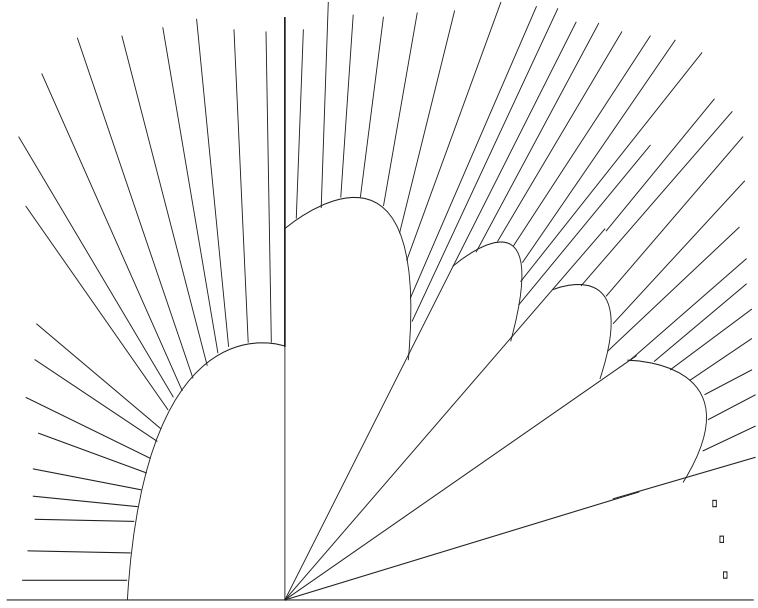} \caption{Comb
Space}

\end{figure}

\section{The Properties Of "COMB" Space}
\begin{itemize}
\item[a)] $dim(\partial X)=0$. For every $n$ we have that $K_n$ is
bounded. That means that we define a finite number of vertices on
every $K_n$ so we add a finite number of geodesic rays. So, all
the infinite geodesic rays are countable. So $\partial X $ is
countable. Now a countable metric space has dimension 0 (see \cite
{H-W} page 18). So dim($\partial X $)=0
\newline
\item[b)] X is a hyperbolic space with the "path" metric. That is
true since every pair of points of $X$ can be joined by a path of
finite length. Also let $l$ be a closed curve of X then $l$ is a
closed curve in $\mathbb H$$^2$ and $length(l)_X \geqslant
length(l)_{{\mathbb H} ^2}$. But since $\mathbb H$$^2$ is
hyperbolic we have the isoperimetric inequality $Area(l) \leqslant
c*length(l)_{{\mathbb H} ^2}$ so $Area(l) \leqslant c*length(l)_X$
which means that X is hyperbolic.(see \cite{GRO1}, \cite{Bow})
\newline
\item[c)] $asdim(X) = 2$. That is because $X$ contains arbitrarily
large balls $B(x,n)\subset \mathbb{H}$$^2$ for every $n \in
\mathbb N$.
\newline
\item[d)] $X$ is a visual hyperbolic space with $D = 1$ since for
every $x$ in $X$ there exists a geodesic from $x_0$ to $x$. Let's
call that $g_1$. If $g_1$ can be extended to infinity then we have
nothing to prove. Let $g_1$ be finite ,then $x$ must belong to a
sector $S(a_n,a_{n+1},N_n)$. We extend $g_1$ until it meets $K_n$
at a point $v_1$. Then by the construction of $X$ there exists an
infinite geodesic $r$ corresponding to the vertex on $K_n$ $v$
such that $d(v_1,v)$ is less than 1. Then obviously $d(x,r)$ is
less than 1.\newline


\end{itemize}

So $X$ is a visual hyperbolic metric space such that that $asdim X
> dim \bd X + 1$.

We remark that it is not very hard to see that $X$ is
quasi-isometrically embedded in $\mathbb{H}$$^2$.

\bibliographystyle{amsplain}

\begin{thebibliography}{99}
\bibitem{GRO1}
M. Gromov, { \em Hyperbolic groups},  Essays in group theory (S. M.
Gersten, ed.), MSRI Publ. 8, Springer-Verlag, 1987 pp. 75-263.
\bibitem {GRO2}
 M.Gromov {\em Asymptotic invariants of
infinite groups }, 'Geometric group theory', (G.Niblo, M.Roller,
Eds.), LMS Lecture Notes, vol. 182, Cambridge Univ. Press (1993)
\bibitem {Bo-S}
M.Bonk and O.Schramm,{\em Embeddings of Gromov Hyperbolic Spaces},
Gafa Geom.Funct.Anal ,Vol 10(2000) ,266-306 .
\bibitem {B-L} S.Buyalo, N.Lebedeva {\em Capacity dimension of
locally self similar spaces}, preprint, August 2005.
\bibitem {H-W}
W.Hurewitz and H.Wallman, {\em'Dimension Theory'} ,Princeton
University Press (1969).
\bibitem {Dr }
A.Dranishnikov{\em Asymptotic topology}, Russian Math.Surveys
55(2000), No 6, 71-116.
\bibitem {Yu}
G.Yu,{\em The Novicov conjecture for groups with finite asymptotic
dimension}, Ann. of Math. 147(1998) ,No 2 , 325-335.
\bibitem {R}
J.Roe, {\em Lectures on Coarse Geometry} AMS University Lecture
Series, 2003
\bibitem{Bow}
B.H.Bowditch {\em A short proof that a Subquadratic Isoperimetric
Inequality Implies a Linear One} , Michigan Math J.42(1995)

\end{thebibliography}

\end{document}